\documentclass[11pt]{article}
\usepackage{amssymb,amsmath,latexsym}
\usepackage[mathscr]{eucal}

\DeclareMathOperator*{\essinf}{ess\,inf}

\usepackage[mathscr]{eucal}
\usepackage{graphicx}
\usepackage{latexsym}
\usepackage{amssymb}
\usepackage{psfrag}
\usepackage{epsfig}

\begin{document}

\newtheorem{lemma}{Lemma}
\newtheorem{theorem}{Theorem}
\newtheorem{corollary}{Corollary}
\newtheorem{definition}{Definition}
\newtheorem{example}{Example}
\newtheorem{proposition}{Proposition}
\newtheorem{condition}{Condition}
\newtheorem{assumption}{Assumption}
\newtheorem{conjecture}{Conjecture}
\newtheorem{problem}{Problem}
\newtheorem{remark}{Remark}

\def\thelemma{\arabic{section}.\arabic{lemma}}
\def\thetheorem{\arabic{section}.\arabic{theorem}}
\def\thecorollary{\arabic{section}.\arabic{corollary}}
\def\thedefinition{\arabic{section}.\arabic{definition}}
\def\theexample{\arabic{section}.\arabic{example}}
\def\theproposition{\arabic{section}.\arabic{proposition}}
\def\thecondition{\arabic{section}.\arabic{condition}}
\def\theassumption{\arabic{section}.\arabic{assumption}}
\def\theconjecture{\arabic{section}.\arabic{conjecture}}
\def\theproblem{\arabic{section}.\arabic{problem}}
\def\theremark{\arabic{section}.\arabic{remark}}

\newcommand{\manualnames}[1]{
\def\thelemma{#1.\arabic{lemma}}
\def\thetheorem{#1.\arabic{theorem}}
\def\thecorollary{#1.\arabic{corollary}}
\def\thedefinition{#1.\arabic{definition}}
\def\theexample{#1.\arabic{example}}
\def\theproposition{#1.\arabic{proposition}}
\def\theassumption{#1.\arabic{assumption}}
\def\theremark{#1.\arabic{remark}}
}

\newcommand{\beginsec}{
\setcounter{lemma}{0} \setcounter{theorem}{0}
\setcounter{corollary}{0} \setcounter{definition}{0}
\setcounter{example}{0} \setcounter{proposition}{0}
\setcounter{condition}{0} \setcounter{assumption}{0}
\setcounter{conjecture}{0} \setcounter{problem}{0}
\setcounter{remark}{0} }

\newcommand{\la}{\lambda}
\newcommand{\eps}{\varepsilon}
\newcommand{\vph}{\varphi}
\newcommand{\al}{\alpha}
\newcommand{\bet}{\beta}
\newcommand{\gam}{\gamma}
\newcommand{\kap}{\kappa}
\newcommand{\s}{\sigma}
\newcommand{\del}{\delta}
\newcommand{\om}{\omega}
\newcommand{\Gam}{\mathnormal{\Gamma}}
\newcommand{\Del}{\mathnormal{\Delta}}
\newcommand{\Th}{\mathnormal{\Theta}}
\newcommand{\La}{\mathnormal{\Lambda}}
\newcommand{\X}{\mathnormal{\Xi}}
\newcommand{\PI}{\mathnormal{\Pi}}
\newcommand{\Sig}{\mathnormal{\Sigma}}
\newcommand{\Ups}{\mathnormal{\Upsilon}}
\newcommand{\Ph}{\mathnormal{\Phi}}
\newcommand{\Ps}{\mathnormal{\Psi}}
\newcommand{\Om}{\mathnormal{\Omega}}

\newcommand{\D}{{\mathbb D}}
\newcommand{\M}{{\mathbb M}}
\newcommand{\N}{{\mathbb N}}
\newcommand{\Q}{{\mathbb Q}}
\newcommand{\R}{{\mathbb R}}
\newcommand{\U}{{\mathbb U}}
\newcommand{\Z}{{\mathbb Z}}

\newcommand{\calA}{{\cal A}}
\newcommand{\calB}{{\cal B}}
\newcommand{\calC}{{\cal C}}
\newcommand{\calD}{{\cal D}}
\newcommand{\calE}{{\cal E}}
\newcommand{\calF}{{\cal F}}
\newcommand{\calG}{{\cal G}}
\newcommand{\calH}{{\cal H}}
\newcommand{\calI}{{\cal I}}
\newcommand{\calJ}{{\cal J}}
\newcommand{\calL}{{\cal L}}
\newcommand{\calM}{{\cal M}}
\newcommand{\calN}{{\cal N}}
\newcommand{\calP}{{\cal P}}
\newcommand{\calR}{{\cal R}}
\newcommand{\calS}{{\cal S}}
\newcommand{\calT}{{\cal T}}
\newcommand{\calU}{{\cal U}}
\newcommand{\calV}{{\cal V}}
\newcommand{\calX}{{\cal X}}
\newcommand{\calY}{{\cal Y}}

\newcommand{\scrA}{\mathscr{A}}
\newcommand{\scrM}{\mathscr{M}}
\newcommand{\scrS}{\mathscr{S}}

\newcommand{\frA}{\mathfrak{A}}
\newcommand{\frS}{\mathfrak{S}}

\newcommand{\bU}{{\bf U}}
\newcommand{\bX}{{\bf X}}
\newcommand{\bZ}{{\bf Z}}

\newcommand{\lan}{\langle}
\newcommand{\ran}{\rangle}
\newcommand{\uu}{\underline}
\newcommand{\oo}{\overline}
\newcommand{\skp}{\vspace{\baselineskip}}
\newcommand{\noi}{\noindent}
\newcommand{\supp}{{\rm supp}}
\newcommand{\diag}{{\rm diag}}
\newcommand{\trace}{{\rm trace}}
\newcommand{\w}{\wedge}
\newcommand{\lt}{\left}
\newcommand{\rt}{\right}
\newcommand{\pl}{\partial}
\newcommand{\abs}[1]{\lvert#1\rvert}
\newcommand{\norm}[1]{\lVert#1\rVert}
\newcommand{\mean}[1]{\langle#1\rangle}
\newcommand{\To}{\Rightarrow}
\newcommand{\til}{\widetilde}
\newcommand{\wh}{\widehat}
\newcommand{\dist}{{\rm dist}}
\newcommand{\grad}{\nabla}

\newcommand{\be}{\begin{equation}}
\newcommand{\ee}{\end{equation}}
\newcommand{\proof}{\noindent{\it Proof.}}
\newcommand{\ink}{\rule{.5\baselineskip}{.55\baselineskip}}
\newcommand{\qed}{\ink}

\newcommand{\mumin}{{\mu_*}}
\newcommand{\rank}{{\rm Rank}}
\newcommand{\Ind}[1]{{\bf 1}_{\left\{ #1 \right\} }}
\def\Half{{\frac 1 2}}

\title{Efficient routing in heavy traffic\\ under
partial sampling of service times}

\author{Rami Atar and Adam Shwartz\\
Department of Electrical Engineering\\
Technion--Israel Institute of Technology\\
\{atar\}\{adam\}@ee.technion.ac.il\\
http://www.ee.technion.ac.il/people/\{adam\}\{atar\}}

\maketitle

\begin{abstract}
We consider a queue with renewal arrivals and $n$ exponential
servers in the Halfin-Whitt heavy traffic regime, where $n$ and the
arrival rate increase without bound, so that a critical loading
condition holds. Server $k$ serves at rate $\mu_k$, and the
empirical distribution of $\{\mu_k\}_{k=1,\ldots,n}$ is assumed to
converge weakly.
We show that very little information on the service rates
is required for a routing mechanism to perform well. More precisely,
we construct a routing mechanism that has access to a {\em single
sample} from the service time distribution of each of
$n^{\Half+\eps}$ randomly selected servers ($\eps>0$), but not to
the actual values of the service rates, the performance of which is
asymptotically as good as the best among mechanisms that have the
complete information $\{\mu_k\}_{k=1,\ldots,n}$.
\\
\\
\noindent {\em Keywords:} Halfin-Whitt regime; routing policies;
service time sampling

\noindent {\em MSC2000:}
Primary: 60F17. Secondary: 68M20, 90B15, 90B22, 60K30, 60K25.
\end{abstract}

\section{Introduction} \label{sec1}
\beginsec

In the many-server parametric regime of Halfin and Whitt
\cite{halwhi}, a critically loaded diffusively scaled system has the
property that the fraction of time when queues are empty is neither
close to 0 nor 1, a situation that is often observed in
applications. Particularly, it has been suggested that this regime
is suitable for modeling large call centers
\cite{gkm},
and various models motivated by this application have been studied, where a
many-server system operates in this regime (see \cite{whi} for a
review). In models that involve heterogenous servers, a principal
problem is to find an efficient routing policy
\cite{arm, ata2,atamanrei, atamansha, bhz,
gurwhi,tez1, daitez}. In all previous works
on routing control in this
regime, the proposed routing mechanisms are assumed to have complete
information about the service rates of each server (where by `rate'
we refer to the parameter of the exponential service time
distribution, assumed by most authors; however see \cite{tez1} for
more general service times). Since often in
applications the routing control mechanism has little knowledge of
the performance of each individual server, it is natural to ask
whether it can perform near optimality with less information on
these parameters. Our goal in this note is to argue that sufficient
information for this purpose is a single sample of service time from
a negligible fraction of the servers.

The pioneering work of Halfin and Whitt \cite{halwhi} considers a
queue with renewal arrivals and identical exponential servers, where
the number of servers and the rate of arrivals are scaled up so that
the queue remains critically loaded. The second order asymptotics of
the process representing the number of customers in the system is
shown to converge to a diffusion. When the servers are heterogenous,
it was shown in \cite{arm} that, in presence of customers of a
single class, the policy which routes jobs to the fastest server
among those that are free at the time of routing (and does not allow
interruption of service) is asymptotically optimal in terms of the
queue length as well as the delay of an arriving customer. Analogous
results are available for the case of random, i.i.d.\ service rates
\cite{ata1} and, under appropriate assumptions, for
hyperexponential service times \cite{tez1}.

We mention that works that characterize the fluid and diffusion
scaling limits are available for homogeneous servers with general
service time distributions \cite{kasram, reed}. The question
that we address here is also very natural in this wider context.
Note, however, that for heterogenous servers with general service
times, an asymptotically optimal routing policy is not known even
when the routing mechanism has access to all service time
distributions (with the exception of \cite{tez1}). For this reason
we confine our treatment to the exponential case.

As mentioned above, we assume that the routing mechanism has access
only to {\it samples} from the service time distribution of some of
the servers. We show that, perhaps counter-intuitively, very little
sampling is required for asymptotically optimal performance: It
suffices to collect a {\it single} sample from each server in a set
of $r$ randomly selected servers, where $r$ is as small as $n^{\Half
+ \eps }$ ($\eps>0$). The proposed policy always routes jobs to
non-sampled servers if such are available, and otherwise, routes to
the server for which the sampled service time is smallest among the
(sampled) servers that are available at the time. It is shown to be
asymptotically optimal in the sense that the diffusion limit of the
process representing the total number of customers in the system
(characterized in Theorem \ref{th1}) is stochastically dominated by
any subsequential limit under any (work conserving, nonanticipating)
policy (see Theorem \ref{th2}). This includes policies that have
access to the complete information on service rates. A similar
statement holds for the queue length processes (simply by
\eqref{04}).

A clear practical advantage of our approach is that it is not
necessary to invest in measuring various characteristics precisely,
or collect accurate information on the performance of the servers.
In addition, the policy proposed has a desired robustness property
in that its performance is nearly optimal regardless of the values
of system parameters, as long as the basic assumptions hold. These
assumptions on the empirical measure of the rates and its first and
second-order limits \eqref{101}--\eqref{117} are quite general, and
so are the assumptions on the limiting distribution.

The proofs are based on an estimate on the number of errors in
ordering the servers according to their sampled data (Lemma
\ref{lem2}), an estimate on the total idle time encountered by
servers that have relatively high priority (Lemma \ref{lem1}), and
the technique developed in \cite{ata1} (proof of Theorem \ref{th1}).

In the next section we describe the model and the proposed policy,
and state the main results. The proofs appear in Section \ref{sec2}.

\section{Model and main results}\label{secres}\beginsec

We fix some notation. Denote by $\D$ the space of functions from
$\R_+$ to $\R$ that are right continuous on $\R_+$ and have finite
left limits on $(0,\infty)$ (RCLL), endowed with the usual Skorohod
topology \cite{bil}. If $X^n$, $n\in\N$ and $X$ are processes with
sample paths in $\D$ (respectively, real-valued random variables) we
write $X^n\Rightarrow X$ to denote weak convergence of the measures
induced by $X^n$ on $\D$ (respectively, on $\R$) to the measure
induced by $X$, as $n\to\infty$. For $X\in\D$ we write
$\abs{X}_{*,t}:=\sup_{0\le s\le t}\abs{X(s)}$. For $x\in\R$, write
$x^+=\max\{x,0\}$ and $x^-=\max\{-x,0\}$.

A complete probability space $(\Om,\calF,P)$ is given, supporting
all random variables and stochastic processes defined below.
Expectation w.r.t.\ $P$ is denoted by $E$.
We consider a single queue fed by renewal arrivals, with parallel exponential
servers. The model is
parameterized by $n\in\N$, where $n$ also represents the number of
servers. The $n$ servers are labeled as $1,\ldots,n$, and, for the
$n$th system, deterministic parameters
$\mu_k^n\in[\uu{\mu},\bar\mu]$ are given, where $\mu^n_k$ represents
service rate of server $k$, and $0<\uu{\mu}\le\bar\mu<\infty$ are
constants independent of $n$. We assume weak convergence of the empirical
measure of $\{\mu^n_k\}$,
 \begin{equation}\label{101}
 L^n=n^{-1}\sum_k\del_{\mu^n_k}\to m,
 \end{equation}
where $m$ is a probability measure on $\R$ (supported on $[\uu{\mu},\bar\mu]$).
The mean is denoted by $\mu=\int xdm$. A second order type
approximation is further assumed on the rate parameters, namely that
the limit
 \begin{equation}\label{102}
 \lim_n n^{-\Half}\sum_{k=1}^n(\mu^n_k-\mu):=\wh\mu
 \end{equation}
exists as a finite number. Denoting $\mu_*=\essinf m$, we finally
assume
\begin{equation}\label{117}
 \lim_{n\to\infty}\#\{k:\mu^n_k<\mu_*-\eps\}n^{-\Half}=0, \quad
 \eps>0.
 \end{equation}

 \begin{example}\label{X1}\rm
A special case of assumptions
\eqref{101}, \eqref{102} and
 \eqref{117} is when there
is a fixed number of pools of servers with $a_in+O(1)$ servers at
pool $i$, and where each server at pool $i$ serves at rate
$b_i+c_in^{-\Half}$ (for constant $a_i,b_i,c_i$; $a_i>0$), a setting
that is common (for example \cite{arm} in a single class setting,
and \cite{daitez} in a multiclass setting).
\end{example}
\begin{example}\rm
We point out that there is more flexibility in the choice of the
parameters. For example, if we have two pools of size $0.2 n +
n^{\frac 3 4} $ and $0.8 n + n^{\frac 4 5} $ with rates $ 1 + 4
n^{-\frac 1 6} + n^{-\Half} $ and, respectively, $ 2 - n^{-\frac 1
6} $, then our assumptions still hold. A more general case is as
follows. We have a fixed number of pools of sizes $ a_i n + f_i (n)
$, with respective rates $ b_i + c_i n^{-\Half } + g_i (n) $. Then
assumptions~\eqref{101}--\eqref{117} hold provided that $f_i (n) = o
(n) $, $g_i (n) = o (1) $ and that the limit
$$
\lim_{n\to\infty}n^{\Half} \sum_i a_i g_i (n)
$$
exists. This is verified by a straightforward, if lengthy
calculation
using $ \mu = \sum_i a_i b_i $.
\end{example}
\begin{example}\rm
It is sometimes very natural to regard the rates $\{\mu_k\}$ as
random variables, and thus to consider the queueing process, as well
as its scaling limit, as processes in random environment. The case
where the service rates are i.i.d.\ random variables, drawn from a
common distribution $m$, was considered in \cite{ata1}. In this
case, the law of large numbers implies that \eqref{101} and
\eqref{117} hold with probability one, and the central limit theorem
implies a variation of \eqref{102}, in which $\wh \mu $ is a normal
random variable. Although we assume throughout that the service
rates are deterministic, we would like to comment that all our
results can be formulated for an i.i.d.\ random environment, with
basically the same proofs.
\end{example}
The initial configuration is now described. Let $Q^n_0$ be a
$\Z_+$-valued random variable, representing the initial number of
customers in the buffer. Let $B_{k,0}^n$, $k=1,\ldots,n$ be
$\{0,1\}$-valued random variables representing the initial state of
each server, where $B_{k,0}^n=1$ if and only if server $k$ initially
serves a customer. We restrict to non-idling policies, so that in
particular $Q^n_0>0$ only if $B_{k,0}^n=1$ for all $k=1,\ldots,n$.
The total number of customers initially in the system
is denoted by $X^n_0=Q^n_0+\sum_{k=1}^nB^n_{k,0}$.
Note that, by assumption, we have the relation $Q_0^n=(X_0^n-n)^+.$
We assume
\begin{equation}\label{23}
\wh X^n_0:=n^{-\Half}(X^n_0-n)\To\xi_0,
\end{equation}
 where $\xi_0$ is a random variable.

To define the arrival process, we are given parameters $\la^n>0$,
$n\in\N$ satisfying $\lim_n\la^n/n=\la>0$, and a sequence of
strictly positive i.i.d.\ random variables $\{\check U(l),l\in\N\}$,
with mean $E\check U(1)=1$ and variance $\check C^2={\rm Var}(\check
U(1))\in[0,\infty)$. With $\sum_1^0=0$, the number of arrivals up to
time $t$ for the $n$th system is given by $
A^n(t)=\sup\{l\ge0:\sum_{i=1}^l\check U(i)/\la^n\le t\}. $ The
arrival rates are further assumed to satisfy the second order
relation
\begin{equation}\label{19}
 \lim_nn^{-\Half}(\la^n-n\la)=\wh\la,
\end{equation}
for some $\wh\la\in\R$. The `heavy traffic' condition on the first
order parameters is assumed, namely
\begin{equation}\label{21}
 \la=\mu,
\end{equation}
indicating that the system is critically loaded.
For each $k=1,\ldots,n$, we let
$B^n_k$ be a stochastic process taking values in $\{0,1\}$,
representing the status of server $k$: when $B^n_k(t)=1$ [resp.,
$0$] we say that server $k$ is busy [resp., idle]. Let
$I^n_k(t)=1-B^n_k(t)$ for $k=1,\ldots,n$, and $t\ge0$. For
$k=1,\ldots,n$, let $R^n_k$ [resp., $D^n_k$] be a $\Z_+$-valued
process with nondecreasing right-continuous sample paths,
representing the number of routings of customers to server $k$
within $[0,t]$ [resp., the number of jobs completed by server $k$ by
time $t$]. Thus
 \begin{equation}\label{05}
 B^n_k(t)=B^n_{k,0}+R^n_k(t)-D^n_k(t),\quad k=1,\ldots,n,\ t\ge0.
 \end{equation}
To describe the processes $D^n_k$, let $\{S_k,k\in\N\}$ be i.i.d.\
standard Poisson processes, each having right-continuous sample
paths. The processes $D^n_k$ are assumed to satisfy
\begin{align}\label{13}
 D^n_k(t) & =S_k(T^n_k(t)),\quad k=1,\ldots,n, \\
\intertext{where}
\label{03}
 T^n_k(t) &=\mu_k^n\int_0^tB^n_k(s)ds,\quad k=1,\ldots,n.
\end{align}
Let $X^n$, $Q^n$ and $I^n$ be defined as
 \begin{align}\label{12}
 X^n(t) =X^n_0+A^n(t)-\sum_{k=1}^nD^n_k(t), \quad
 Q^n(t) =Q^n_0+A^n(t)-\sum_{k=1}^nR^n_k(t), \quad
 I^n(t) =\sum_{k=1}^nI^n_k(t).
\end{align}
These processes represent the number of customers in the system, the
number of customers in the buffer and, respectively, the number of
servers that are idle.

The routing policy, that will be described below, does not have
access to the service rates $\mu^n_k$, but it has access to samples
from the service time of $r$ of the servers, selected at random, and
no information at all on service rates of the others. More
precisely, let $r=r_n\in\N$, $r\le n$ be given, and let
$\Sig=\Sig^n$ be a random variable uniformly distributed over the
set of all subsets of $\{1,\ldots,n\}$ that have cardinality $r$. We
denote $\Sig^c=\{1,\ldots,n\}\setminus\Sig$. For each $k\in \Sig$,
let $\s_k=\s_k^n$ be an independent copy drawn from the service time
distribution of server $k$. That is, $\s_k$ is an exponential random
variable with parameter $\mu^n_k$ and, conditioned on $\Sig$,
$\{\s_k\}_{k\in\Sig}$ are independent. We choose
\begin{equation}\label{103}
 r_n=[n^{\beta_r}],
\end{equation}
where $\beta_r\in(\Half,1]$. Denote $\wh\mu_k=1/\s_k$, $k\in\Sig$.

The four stochastic primitives introduced, as listed below, are
assumed to be mutually independent, for each $n$:
 \begin{equation}\label{104}
 \lt(X_0^n,
 \{B^n_{k,0}\}_{k=1,\ldots,n}\rt),\quad
 \{S_k\}_{k\in\N},\quad A^n,
 \quad \Big(\Sig^n,\{\s_k^n\}_{k\in\Sig^n}\Big).
 \end{equation}

Routing is based on an ordering of the servers according to whether
they are in $\Sig$ and, within $\Sig$, according to the value of
$\wh\mu_k$. A permutation $\rank=\rank_n$ of $\{1,\ldots,n\}$ is
defined as follows.
On the probability-one event that
the $\wh\mu_k$ are all distinct, the
set $\Sig$ is mapped by $\rank$ onto $\{1,\ldots,r\}$
(and $\Sig^c$ onto $\{r+1,\ldots,n\}$). For $k,l\in\Sig$,
$\rank(k)<\rank(l)$ if and only if $\wh\mu_k<\wh\mu_l$. For
$k,l\in\Sig^c$, $\rank(k)<\rank(l)$ if and only if $k<l$.

The routing policy favors servers ranked higher (namely those that
have high value under the map $\rank$). That is, when a customer
arrives to the system to find
more than one idle server,
the customer is routed to the server with
highest rank among those servers.
Since it is assumed that the routing policy is work
conserving (non-idling),
when the queue is nonempty and a server has just
finished serving, a customer (from the head of the line) is routed
to this server, and
when a customer arrives to the system to find
exactly one server that is idle, it is
instantaneously routed to that server. As a result,
 \begin{equation}\label{04}
 Q^n(t)=(X^n(t)-n)^+,\quad I^n(t)=(X^n(t)-n)^-
 \end{equation}
 holds for all $t$.
Also, service is non-interruptible, in the sense that a customer
completes service at the server it is first assigned.

This completes the description of the process
$$\PI_0^n:=(\{B_k^n\},\{R^n_k\},\{D^n_k\},X^n,Q^n,I^n).$$ It can be
seen that this description uniquely determines $\PI_0^n$. We
sometimes refer to this process as {\it policy $\PI_0^n$}. Later we
use some of the symbols above (such as $X^n$) to denote quantities
that have the same meaning (such as the number of customers in the
$n$th system) under a different routing policy $\PI^n$. To avoid confusion,
we therefore make specific reference to policy $\PI_0^n$ when
necessary.

Finally, we make a simplifying assumption about the initial
occupation of servers, namely that only servers that are ranked low
may initially be idle:
 \begin{equation}
   \label{109}
   B^n_{k,0}=1_{\{\rank(k)>I^n_0\}},
 \end{equation}
where
 \begin{equation}\label{110}
 I^n_0=(X^n_0-n)^-
 \end{equation}
is the initial number of idle servers.

Let $\wh X^n$ be a centered, normalized version of the process
$X^n$, defined by
 \begin{equation}\label{113}
 \wh X^n=n^{-\Half}(X^n-n).
 \end{equation}
Our main result is the following.
\begin{theorem}\label{th1}
Under policy $\PI_0^n$, the processes $\wh X^n$ converge weakly to
the unique solution $\xi$ of
 \begin{equation}\label{85}
 \xi(t)=\xi_0+\s w(t)+\beta t+\mumin\int_0^t\xi(s)^-ds,\quad t\ge0,
 \end{equation}
where $\s^2=\mu(\check C^2+1)$, $\beta=\wh\la-\wh\mu$, and $w$ is a
standard Brownian motion, independent of $\xi_0$.
\end{theorem}

The result above is to be compared with Proposition 4.2 of
\cite{arm} and Theorem 2.2 of \cite{ata1} (for the case of a finite
number of server pools and, respectively, random environment). In
these references, equation \eqref{85} arises in the limit under a
policy defined similarly to $\PI_0$, but where the servers are
ordered according to the actual values of $\mu_k$, $k=1,\ldots,n$.
In \cite{arm} it is further shown that this policy asymptotically
achieves the best performance in a large class of routing policies.
Because our setting is different than \cite{arm}, we will state and
prove an analogous result, so as to exhibit that $\PI_0$ is
asymptotically optimal.

Toward this end, let us first comment on an alternative
representation of the departure process. By \eqref{13}, this process
is given as $\sum_{k=1}^nD^n_k(t)=\sum_{k=1}^nS_k(T^n_k(t))$, where
$S_k$ are independent standard Poisson processes. In fact, the
departure process can also be represented as
 \begin{equation}\label{112}
 \sum_{k=1}^nD^n_k(t)=S^n\Big(\sum_{k=1}^nT^n_k(t)\Big),
 \end{equation}
where, for every $n$, $S^n$ is a standard Poisson process,
independent of the remaining primitive data, that is, of the first,
third and fourth items of \eqref{104}. This statement (along with a
variation of it, stated in Section \ref{sec2}) is due to a standard
superposition argument for Poisson processes, for which the reader is
referred to Proposition 3.1 of \cite{ata1}.

We now define a class of policies by keeping the description of this
section but abandoning the specifics of the routing mechanism. More
precisely, we write $\PI^n\in\calP^n$ for any process
$$
 \PI^n=(\{B_k^n\},\{R^n_k\},\{D^n_k\},X^n,Q^n,I^n)
$$
satisfying all relations stated throughout this section, from its
beginning to the statement of Theorem \ref{th1}, save the two
paragraphs following display \eqref{104}, and satisfying, in
addition, work conservation \eqref{04} and the representation
\eqref{112}, for some standard Poisson processes $S^n$, independent
of the remaining primitive data. Note, in particular, that the
routing mechanism may have access to $\{\mu_k\}$.
See Remark \ref{rem1} about the role played by the work conservation
condition \eqref{04}.
 We refer to any element of $\calP^n$ as a {\it policy}.

\begin{theorem}
  \label{th2}
For $n\in\N$ and any policy $\PI^n\in\calP^n$, let $\wh X^n$
be the normalized version \eqref{113} of the corresponding process
$X^n$. Then
there exist processes $\X^n$ that converge weakly, as
$n\to\infty$, to the solution $\xi$ to \eqref{85} and
$$\wh X^n(t)\ge\X^n(t),\quad t\ge0, \text{ $P$-a.s.},\,
n\in\N.$$
\end{theorem}

Since by Theorem \ref{th1}, $\xi$ is obtained at the limit under
$\PI^n_0$, the result above demonstrates that $\PI^n_0$
asymptotically optimal.

\section{Proofs}\label{sec2}
\beginsec

We begin with the following.

\begin{lemma}\label{lem2}
Let $0<\phi<\psi<\infty$, $\beta>\Half$, $c_1>0$ and $c_2>0$ be
given constants. For $n\in\N$ denote $\ell^1=\ell^1_n=[c_1n^\beta]$,
$\ell^2=\ell^2_n=[c_2n^\beta]$, and let
$\phi^n_1,\ldots,\phi^n_\ell$ and $\psi^n_1,\ldots,\psi^n_\ell$ be
positive real numbers with
 $$
 \sup_{n,i}\phi^n_i\le\phi<\psi\le\inf_{n,i}\psi^n_i.
 $$
For $n\in\N$ and $i\in\{1,\ldots,\ell\}$ let $\Ph^n_i$ [resp.,
$\Ps^n_i$] be an exponential random variable with parameter
$\phi^n_i$ [resp., $\psi^n_i$]. For each $n$ assume that
$\{\Ph^n_i\}$ are mutually independent and that so are
$\{\Ps^n_i\}$. Then there exist $\gamma>0$ and $\kap>0$ such that,
with $ \theta_n = \gamma\log n $, one has
\begin{equation}\label{121}
\lim_{n\to\infty}P \left( \sum_{i=1}^{\ell_n^1} \Ind{\Ph^n_i \geq
\theta_n } \leq n^{\Half + \kap }
  \right ) =0,
\end{equation}
\begin{equation}\label{122}
\lim_{n\to\infty}P \left ( \sum_{i=1}^{\ell_n^2} \Ind{\Ps^n_i \geq
\theta_n } \geq n^{\Half - \kap }
  \right ) = 0 .
\end{equation}
\end{lemma}

\noindent \proof {} By stochastic monotonicity of the exponential
random variable with respect to its parameter, it clearly suffices
to prove the claim for the case $\phi^n_i=\phi$, $\psi^n_i=\psi$,
all $n$ and $i$. To prove the claim under this assumption, let $\kap
$ and $\gamma$ be strictly positive constants satisfying
 \begin{equation}\label{123}
 \phi\gamma<\beta-\Half-\kap<\beta-\Half+\kap<\psi\gamma.
 \end{equation}
Write $\Ph_i$ [resp. $\Ps_i$] for $\Ph^n_i$ [resp., $\Ps^n_i$]. Then
for any $\al>0$ and $\{\theta_n\}$,
\begin{align*}
P \left ( \sum_{i=1}^{\ell_n^1} \Ind{\Ph_i \geq \theta_n } \leq
n^{\Half + \kap }
  \right )
& \leq e^{\alpha n^{\Half + \kap }}
     E \exp \left [ -\alpha \sum_{i=1}^{\ell_n^1} \Ind{\Ph_i \geq \theta_n }
            \right ]         \\
& =  e^{\alpha n^{\Half + \kap }}
     \left ( P (x_1 < \theta_n ) + e^{-\alpha } P (x_1 \geq \theta_n )
     \right )^{\ell_n^1}                                                    \\
& =  e^{\alpha n^{\Half + \kap }}
     \left ( 1 - e^{-\phi\theta_n} + e^{-\alpha } e^{-\phi\theta_n} \right)^{\ell_n^1} \\
& =  e^{\alpha n^{\Half + \kap }}
     \exp \left\{ \ell_n^1 \log \left ( 1 + e^{-\phi\theta_n} ( e^{-\alpha } - 1 )
           \right ) \right \}                                              \\
& \leq \exp \left\{ \alpha n^{\Half + \kap } -
                    \ell_n^1 \left ( e^{-\phi\theta_n} (1- e^{-\alpha } )
           \right \} \right )  \ .
\end{align*}
For the last expression to converge to zero, we need
 $ K_n:=\alpha
n^{\Half + \kap } -
      \ell_n^1 \left ( e^{-\phi\theta_n} (1- e^{-\alpha }) \right ) \to -\infty .$
Fix $\nu>0$ and set $ \theta_n = \gamma \log n$, $\al=\alpha_n = \nu
\log n . $ Then
 $$
K_n
 \leq \nu \log n \cdot n^{\Half + \kap } -
     (c_1 n^\beta -1) \left (  n^{- \phi\gamma } ( 1 - n^{ - \nu })
     \right).
 $$
Since by \eqref{123} $ \Half + \kap < \beta - \phi\gamma $, we have
$K_n\to-\infty$, as desired, and thus \eqref{121} holds.

Fix $\eta > 0 $. Since $ (1+x)^k \leq e^{kx} $ for $ x
> -1 $, we have
\begin{align*}
P \left ( \sum_{i=1}^{\ell_n^2} \Ind{\Ps_i \geq \theta_n } \geq
n^{\Half - \kap } \right )
& \leq e^{-\eta n^{\Half - \kap }}
       E \exp\left[ {\eta \sum_{i=1}^{\ell_n^2} \Ind{\Ps_i \geq \theta_n } }
\right]        \\
& =    e^{-\eta n^{\Half - \kap }}
       \left ( 1 - e^{-\psi \theta_n} + e^{-\eta } e^{-\psi \theta_n}
       \right)^{\ell_n^2} \\
& \leq e^{-\eta n^{\Half - \kap }}
       \exp \left\{ \ell_n^2 \left ( e^{-\psi \theta_n}
       ( e^{-\eta } - 1 ) \right ) \right \} \\
&= \exp \left \{ -\eta n^{\Half - \kap }  + c_2 n^\beta n^{-\psi
\gamma } ( e^{-\eta } - 1 ) \right \},
 \end{align*}
where on the last line above we substituted $ \theta_n = \gamma \log
n $. The expression on the last line converges to zero because $
\Half - \kap
> \beta - \psi \gamma$ by \eqref{123}, and \eqref{122} follows.
 \qed

\begin{remark}\rm
 {\bf (a)}
The convergence in \eqref{121}, \eqref{122} is
at a geometric rate, as the proof shows. Thus, by the
Borel-Cantelli lemma, both events
occur for only a finite number of $n$, with probability one.
\\
{\bf (b)} As can be seen in the proof, $\kap $ and $\gamma $ depend
only on $\beta ,\phi $ and $\psi $ (cf.\ \eqref{123}), and not on
$\{c_i\}$.
\end{remark}

Fix $\eps>0$ and let $\al\in(\mumin,\mumin+\eps)$ be a continuity
point of $x\mapsto m([0,x])$.
In what follows, the symbols $n$ and $\eps $ are
omitted from the notation of all
random variables and stochastic processes, and from the parameters
$\mu_k^n$.
Let $M_0=[\uu{\mu},\mu_*-\eps)$,
$M_1=[\mu_*-\eps,\al)$ and $M_2=[\al,\bar\mu]$ (where $[a,b)$ and
$[a,b]$ are interpreted as the empty set if $a>b$), and set
$$
 K_i=\{k\in\{1,\ldots,n\}:\mu_k\in M_i\},\quad i=0,1,2.
$$
Denote
 \begin{equation}\label{75}
I^{(i)}(t)=\sum_{k\in K_i}I_k(t),\quad
T^{(i)}(t)=\sum_{k\in K_i}T_k(t),\quad i=0,1,2.
 \end{equation}
Let also $\wh I^{(i)}=n^{-\Half}I^{(i)}$. By \eqref{13} and
\eqref{12},
 \begin{equation}\label{111}
\wh X(t) =
 \wh X_0+n^{-\Half}A(t)-n^{-\Half}\sum_{k=1}^n S_k(T_k(t)).
 \end{equation}
By a superposition argument for Poisson processes (cf.\ Proposition 3.1
of \cite{ata1}),
\begin{align}
\wh X(t) & = \label{36}
 \wh X_0+n^{-\Half}A(t)-n^{-\Half}\sum_{i=0}^2 S^{(i)}(T^{(i)}(t)),
\end{align}
where $S^{(i)}$, $i=0,1,2$ are standard Poisson processes, mutually
independent, and independent of the first, third and fourth items of
\eqref{104}. In particular,
 \begin{equation}\label{74}
D^{(i)}(t):=\sum_{k\in K_i}D_k(t)=S^{(i)}(T^{(i)}(t)),\quad i=0,1,2.
 \end{equation}
A calculation based on \eqref{21} and \eqref{36} shows (see a
detailed derivation at the end of this section)
 \begin{align}
 \wh X(t)=\wh X_0+W(t)+bt+F(t),\label{30}
\end{align}
where we recall that all quantities depend on $n$ and $\eps $, and
where
 \begin{align}\label{31}
 W(t) &=
         \wh A(t)-\sum_{i=0}^2W^{(i)}(t),        \\
\label{47}
 \wh A(t) &=n^{-\Half}(A(t)-\la^nt),        \\
\label{52}
 W^{(i)}(t) &=n^{-\Half}(S^{(i)}(T^{(i)}(t))-T^{(i)}(t)),\quad
 i=0,1,2,        \\
\label{32}
b &=n^{-\Half}(\la^n-n\la)-n^{-\Half}\sum_{k=1}^n(\mu_k-\mu),        \\
\label{33}
F(t) &=
       n^{-\Half}\int_0^t\sum_{k=1}^n\mu_kI_k(s)ds.
 \end{align}

\begin{lemma}
  \label{lem1}
Under $\PI_0^n $, given $\bar t>0$ and $\eps>0$,
 \begin{equation}\label{78}
 |\wh I^{(2)}|_{*,\bar t}\to0 \text{ in probability, as }
 n\to\infty.
 \end{equation}
\end{lemma}

\proof{}
 Step 1:
We will show here that there is a (deterministic) sequence $a_n$
increasing to infinity, so that $a_nn^{\Half}\le r_n$, and such
that, out of the $a_nn^{\Half}$ servers ranked lowest, the number of
those that are in $K_2$ is $o(n^{\Half})$, in the following sense:
 \begin{equation}\label{105}
 \frac{\#\{k\in K_2:\rank(k)\le a_nn^{\Half}\}}{n^{\Half}}\To0.
 \end{equation}

We will apply Lemma \ref{lem2}. To this end let $\phi\in(\mu_*,\al)$
be a continuity point of $x\mapsto m([\uu\mu,x])$. Let $\psi=\al$.
Let $\til K=\{k\in\{1,\ldots,n\}:\mu_k\le\phi\}$. Since
$m([\mu_*,\phi])>0$, it follows from \eqref{101} that,
for some constant $c>0$ and with probability increasing to $1$,
the cardinality of $\til K$ is at least $cn$.
Since the subset $\Sig$ is uniformly
distributed and the number of samples satisfies \eqref{103}, it
follows that, on some events $\Om^n$ satisfying $P(\Om^n)\to1$,
 $$
 \#\til K\cap\Sig\ge c_1n^{\beta_r},\quad
 \# K_2\cap\Sig\le\#\Sig=n^{\beta_r},
 $$
 for a constant $c_1>0$.
 Recall $\wh\mu_k=1/\s_k$, the reciprocal to the sampled service
time.
 We apply Lemma \ref{lem2} with $\Ph_i$ being the samples $\s_k$ with
 index set $\til K\cap\Sig$, and $\Ps_i$ those with index set
 $K_2\cap\Sig$. We obtain, that on an event $\Om_1^n\subset\Om^n$, which also
 satisfies $\lim_{n\to\infty}P(\Om^n_1)=1$,
 \begin{align}\label{124}
 \#\{k\in\til K\cap\Sig:\wh\mu_k\le1/\theta_n\}>
 n^{\Half+\kap},         \\
\label{125}
 \#\{k\in K_2\cap\Sig:\wh\mu_k\le1/\theta_n\}<
 n^{\Half-\kap},
 \end{align}
 where, without loss of generality, $\Half<\Half+\kappa<\beta_r$.
Now, \eqref{124} and the way the map $\rank$ is defined, imply that
all servers $k$ with $\rank(k)\le n^{\Half+\kap}$ have
$\wh\mu_k\le1/\theta_n$ and are in $\Sig$. As a result,
 \begin{align*}
 \#\{k\in K_2:\rank(k)\le n^{\Half+\kap}\} &=\#\{k\in K_2\cap\Sig:\rank(k)\le
 n^{\Half+\kap}\}\\
 &\le \#\{k\in K_2\cap\Sig:\wh\mu_k\le1/\theta_n\}\\
 &\le n^{\Half-\kap},
 \end{align*}
by \eqref{125}. This proves \eqref{105} with $a_n = n^\kap $.

Step 2:
Denote $$K'=\{k\in\{1,\ldots,n\}:\rank(k)>a_nn^{\Half}\},$$ and
$I'=\sum_{k\in K'}I_k$, $T'=\sum_{k\in K'}T_k$, $D'=\sum_{k\in
K'}D_k$. An argument as the one following equation \eqref{111} shows
that $S'(T'(t))=D'(t)$, $t\ge0$, where $S'$ is a standard Poisson
process. Set $\wh S'(t)=n^{-\Half}(S'(nt)-nt)$ and $\wh
I'=n^{-\Half}I'$. We shall show that
 \begin{equation}
   \label{108}
  |\wh I'|_{*,\bar t}\to0 \text{ in probability, as }
  n\to\infty.
 \end{equation}

Note first that the probability of the event $\eta_1:=\{I'(0)=0\}$
converges to one as $n\to\infty$. Indeed, by \eqref{109},
$B_{k,0}=1$ for all $k$ with $\rank(k)>I_0$. By \eqref{23} and
\eqref{110}, $I_0<a_nn^{\Half}$ with probability converging to 1 as
$n\to\infty$. Thus, with probability converging to 1, all servers
$k\in K'$ are initially busy, namely $P(\eta_1)\to1$ as
$n\to\infty$.
Let
 \begin{equation}\label{22}
 \wh S(t)=n^{-\Half}(S(nt)-nt),\quad t\ge0,
 \end{equation}
 where $S$ is a standard Poisson process. It is well known
(cf.\ Lemmas 2 and 4(i) of \cite{atamanrei}) that both $\wh A$ (of
\eqref{47}) and $\wh S$ converge weakly to a zero mean Brownian
motion with diffusion coefficient $\la^{\Half}\check C$, and
respectively, 1.

Given $\gam>0$, consider the event $\eta:=\{|I'|_{*,\bar t}>2\gam
n^{\Half}\}$. On the event $\eta\cap\eta_1$ one can find $0\le
s<t\le\bar t$ such that $I'(y)>0$ for all $y\in[s,t]$, and
$I'(t)-I'(s)>\gam n^{\Half}$. Since the servers in $K'$ are all
ranked higher than those in the complement set, the routing policy
assigns all arrivals within $[s,t]$ to $K'$ servers. Hence by
\eqref{05}, \eqref{13}
and using $ B_k = 1 - I_k $, we have
$$
 \gam n^{\Half}<I'(t)-I'(s)=D'(t)-D'(s)-A(t)+A(s),
$$
and therefore
\begin{multline*}
 \gam <\wh S'(n^{-1}T'(t))-\wh
 S'(n^{-1}T'(s))-\wh A(t)+\wh A(s)\\
 +\sum_{k\in
 K'}\mu_k\int_s^t \wh B_k(y)dy-\la
 n^{\Half}(t-s)-n^{-\Half}(\la^n-\la)(t-s).
\end{multline*}
We have by \eqref{03} and \eqref{75} that $n^{-1}T^{(2)}(t)\le
\bar\mu\bar t=:\tau$. Also, by \eqref{19}, the last term above is
bounded by $c(t-s)$ for some constant $c$ independent of $n$ and
$\eps$.  Let
$$
 \bar w_{\tau}(x,z)=\sup_{|s-t|\le z;s,t\in[0,\tau]}|x(s)-x(t)|,\quad z>0,
$$
denotes the modulus of continuity for $x:[0,\tau]\to\R$. Define
$C(n,\eps)=n^{-1}\sum_{k\in K'}\mu_k-\la$. Then on the event
$\eta\cap\eta_1$, with $\del=t-s$, we have
 \begin{equation}\label{90}
\gam<\bar w_\tau(\wh S',2\bar\mu\del)+\bar w_{\bar t}(\wh
A,\del)+n^{\Half}C(n,\eps)\del+c\del.
 \end{equation}
By \eqref{101}, \eqref{102}, \eqref{21} and the definition of $K'$,
\begin{align*}
 n^{\Half}C(n,\eps)& \le c_1-n^{-\Half}\sum_{k:\rank(k)\le n^{\Half}a_n}\mu_k
 \le c_1-\uu\mu a_n\le-c_2a_n,
\end{align*}
for constants $c_1$, $c_2>0$ and sufficiently large $n$. Hence
$$
 P(|\wh I'|_{*,\bar t}>2\gam)=P(\eta)\le p_1(n,\eps,\gam)+p_2(n,\eps,\gam)+P(\eta_1^c),
$$
where
$$
 p_1(n,\eps,\gam)=P(\text{there exists
 }\del\in(0,a_n^{-\Half}) \text{ such that \eqref{90} holds}),
$$
$$
 p_2(n,\eps,\gam)=P(\text{there exists
 }\del\in[a_n^{-\Half},\bar t] \text{ such that \eqref{90} holds}).
$$
Note that
$$
p_1(n,\eps,\gam)\le P(\bar w_\tau(\wh S',2\bar\mu a_n^{-\Half})+\bar
w_{\bar t}(\wh A,a_n^{-\Half})\ge\gam/2),
$$
$$
 p_2(n,\eps,\gam)\le P(\bar w_\tau(\wh S',2\bar\mu\bar t)+\bar w_{\bar t}(\wh
A,\bar t)\ge-c\bar t+c_2a_n^{\Half}).
$$
Since $\wh S'$ and $\wh A$ converge to processes with continuous
sample paths, both expressions converge to zero as $n\to\infty$.
Since $\lim_nP(\eta_1^c)=0$ and $\gam>0$ is arbitrary, \eqref{108}
follows.

Step 3:
Since $K_2\subset K'\cup((K')^c\cap K_2)$, we have
$$
 \wh I^{(2)}(t)=\frac{I^{(2)}(t)}{n^{\Half}}\le \frac{I'(t)}{n^{\Half}}+\frac{[\#(K'\cup((K')^c\cap K_2))]}{n^{\Half}}\,\bar t,\quad t\in[0,\bar t].
$$
By Step 1 (display \eqref{105}), the last term on the above display
converges to zero in probability. Thus by Step 2 (display
\eqref{108}), statement \eqref{78} follows. This completes the proof
of the lemma.
 \qed

\skp

\noi{\it Proof of Theorem \ref{th1}.} Based on Lemmas \ref{lem2} and
\ref{lem1}, the proof is similar to that of Theorem 2.2 of
\cite{ata1} (only slightly simpler). We include it for completeness
and because the proof of Theorem \ref{th2} is based on it. By
\eqref{30} and \eqref{33}, one has
 \begin{equation}\label{79}
 \wh X(t)=\wh X_0+W(t)+bt+\mumin\int_0^t\wh X(s)^-ds+e(t),
 \end{equation}
 (where all the above quantities depend on $n$) and, with $\wh I_k=n^{-\Half}I_k$,
 \begin{equation}\label{80}
 e(t)=\sum_{k=1}^n(\mu_k-\mumin)\int_0^t\wh I_k(s)ds.
 \end{equation}
 Fix $\bar t>0$. By \eqref{102}, \eqref{19} and \eqref{32},
$b\to\beta=\wh\la-\wh\mu$. We show that the random variables
$\{|W^{(i)}|_{*,\bar t},i=0,1,2,\,n\in\N\}$ are tight. By \eqref{75}
and \eqref{03}, for $i=0,1,2$,
\begin{align}\label{35}
 n^{-1}T^{(i)}(t) &=n^{-1}\sum_{k\in K_i}\mu_kt
 -n^{-1}\sum_{k\in K_i}\mu_k\int_0^tI_k(s)ds.
\end{align}
Hence $0\le n^{-1}T^{(i)}(t)\le\bar\mu\bar t$ for $t\le\bar t$ and
all $n$. Thus by \eqref{52}, $|W^{(i)}|_{*,\bar t}\le|\wh
S^{(i)}|_{*,\bar\mu\bar t},$ where $\wh
S^{(i)}(t)=n^{-\Half}(S^{(i)}(nt)-nt)$. Recall from the proof of
Lemma \ref{lem1} that $\wh S^{(i)}$ converge to a Brownian motion.
Hence $|W^{(i)}|_{*,\bar t}$ are tight.

Next, note that $|e(t)|\le\bar\mu\int_0^t|\wh X(s)|ds.$ Thus the
boundedness of $b$, the tightness of the random variables $\wh X_0$,
$|W^{(i)}|_{*,\bar t}$ and $|\wh A|_{*,\bar t}$, $n\in\N$ (as
follows from the convergence of $\wh A$), and an application of
Gronwall's lemma on \eqref{79}, by which
$
 |\wh X|_{*,\bar t}\le (|\wh X_0|+|W|_{*,\bar t}+|b|\bar
t)\exp(2\bar\mu\bar t),
$
imply that $\{|\wh X|_{*,\bar t},n\in\N\}$ are tight. Since by
\eqref{04}, $\wh I=\wh X^-$, we have that $\{|\wh I|_{*,\bar
t},n\in\N\}$ are tight.

The supremum over $t\le\bar t$ of the absolute value of the last
term in \eqref{35} converges to zero in probability, since $\mu_k$
are assumed to be bounded and $|\wh I|_{*,\bar t}$ are tight. Also,
since $\al$ is a continuity point of $x\mapsto m([0,x])$, we have
that
$$
 n^{-1}\sum_{k\in K_i}\mu_k\to \int_{M_i}xdm=:\rho_i,\quad i=0,1,2.
$$
Note that $\rho_0=0$. As a result,
$n^{-1}(T^{(0)},T^{(1)},T^{(2)})\to\til\rho$ in probability,
uniformly on $[0,\bar t]$, where $\til\rho(t)=(0,\rho_1t,\rho_2t)$.
Recall that $(\wh A,\wh S^{(0)},\wh S^{(1)},\wh S^{(2)})$ are
mutually independent, and that $\wh S^{(i)}$ [resp., $\wh A$]
converges to a standard Brownian motion [a zero mean Brownian motion
with diffusion coefficient $\la^{\Half}\check C$] (see comment
following \eqref{22}). Thus \eqref{31}, \eqref{52} and the lemma on
random change of time \cite[p.\ 151]{bil} show that $W$ converges
weakly to $\s w$, in the uniform topology on $[0,\bar t]$, where $w$
is a standard Brownian motion and $\s^2=\la\check
C^2+\rho_1+\rho_2=\la\check C^2+\mu=\mu(\check C^2+1)$.

By the Skorohod representation theorem, we can assume without loss
of generality that the random variables $\wh X_0$ and $\xi_0$ and
the processes $W$ and $w$ are realized in such a way that, $P$-a.s.,
 \begin{equation}\label{42}
(\wh X_0,W)\to(\xi_0,\s w)\quad \text{ as } n\to\infty.
 \end{equation}
Let $\xi$ be the unique strong solution to equation \eqref{85}. Then
by \eqref{85}, \eqref{79}, the inequality $|x^--y^-|\le|x-y|$, and
Gronwall's inequality,
 \begin{equation}\label{107}
|\wh X-\xi|_{*,\bar t}\le (|\wh X_0-\xi_0|+|b-\beta|+|W-\s
w|_{*,\bar t}+|e|_{*,\bar t})\exp(\mumin\bar t).
 \end{equation}
Now, by \eqref{80}, for $n$ sufficiently large,
 \begin{equation}\label{106}
 |e|_{*,\bar t}\le\eps\bar t\,|\wh
 I|_{*,\bar t}+\bar\mu\bar t\,|\wh I^{(2)}|_{*,\bar t}+(\mu_*-\uu{\mu})\bar t\,|\wh I^{(0)}|_{*,\bar t}.
 \end{equation}
By \eqref{117}, the last term above converges weakly to 0. Combining
\eqref{78}, \eqref{42}, \eqref{107} and \eqref{106},
$$
 \limsup_n P(|\wh X-\xi|_{*,\bar t}>\eps^{\Half})\le\limsup_nP(c\eps|\wh
 I|_{*,\bar t}>\eps^{\Half}),
$$
where $c\in(0,\infty)$ is a constant independent of $n$ and $\eps$.
Note that the law of $|\wh I|_{*,\bar t}$ does not depend on $\eps$.
Hence by tightness of $\{|\wh I|_{*,\bar t},n\in\N\}$ the r.h.s.\ in
the above display converges to zero as $\eps\to0$. Thus $|\wh
X-\xi|_{*,\bar t}\to0$ in probability. Since $\bar t$ is arbitrary,
we have $\wh X\To\xi$.
 \qed

\skp

\noi{\it Proof of Theorem \ref{th2}.} By \eqref{117} there exists a
sequence $\del_n>0$ tending to zero such that
$\zeta_n:=\#\{k:\mu^n_k<\mu_*-\del_n\}n^{-\Half}\to0$. Note that
\eqref{79}, \eqref{80} still hold. Define $\X_1$ as the solution to
 \begin{equation}\label{114}
 \X_1(t)=\wh X_0+W(t)+bt+\mumin\int_0^t\X_1(s)^-ds.
 \end{equation}
Then by \eqref{79}, $\Del:=\wh X-\X_1$ is differentiable, and, using
the inequality $a^--b^-\le-(a-b)^+$ for $a,b\in\R$, we have
$\Del(0)=0$, and
 $$
 \dfrac{d}{dt}\Del(t)\ge-\mumin\Del(t)^++\frac{d}{dt}e(t).
$$
Since $\wh I_k\le n^{-\Half}$ for each $k$, we have by \eqref{80}
$$
 \dfrac{d}{dt}e(t)\ge-v_n,\qquad v_n:=\del_n|\wh I|_{*,\bar t}+\mu_*\zeta_n,
 $$ and $\Del(0)=0$.
By comparison with the ordinary differential equation $du/dt=-\mu_*
u^+-v_n$, $u(0)=0$, we obtain that $\Del(t)\ge-v_nt$, $t\le \bar t$.
Hence $X(t)\ge\X(t)$, where we define $\X(t)=\X_1(t)-v_nt$,
$t\le\bar t$.

It thus remains to show that $\X\To\xi$. For this let us review the
proof of Theorem \ref{th1}. Rather than three processes $D^{(i)}$
\eqref{74} and correspondingly $W^{(i)}$, $i=0,1,2$, \eqref{52}, we
now have a single process $D=\sum_kD_k$ given in terms of a single
standard Poisson process $S^n$ (cf.\ \eqref{112}). The adaptation of
relation \eqref{52} to a single process $W$ is obvious. The
arguments in the proof of Theorem \ref{th1} that lead to the
tightness of $|\wh I|_{*,\bar t}$ and the convergence of $W$ to $\s
w$ hold with obvious modifications. As in that proof, we deduce that
\eqref{42} can be assumed without loss of generality. Equations
\eqref{85}, \eqref{114} and Gronwall's inequality thus yield
$$
|\X_1-\xi|_{*,\bar t}\le (|\wh X_0-\xi_0|+|b-\beta|+|W-\s w|_{*,\bar
t})\exp(\mumin\bar t).
$$
Hence \eqref{42} and the convergence of $b$ to $\beta$ imply that
$\X_1$ converges in probability to $\xi$ uniformly over $[0,\bar
t]$. The random variables $v_n$ converge to zero by tightness of
$|\wh I|_{*,\bar t}$. Since $\bar t$ is arbitrary, we thus obtain
that $\X\To\xi$. This completes the proof of the theorem.
 \qed

\begin{remark}\label{rem1}\rm
Note that the non-idling property is used in the proof of Theorem
\ref{th1} (on which the above proof is based) for deducing tightness
of $|\wh I|_{*,\bar t}$ from that of $|\wh X|_{*,\bar t}$. As can be
easily seen, $|\wh I|_{*,\bar t}$ are not in general tight if the
restriction to non-idling policies is removed.
\end{remark}

\noi{\it Derivation of equation \eqref{30}.}
Recall that $n$ is omitted from some of the notation.
By \eqref{36}, \eqref{47}--\eqref{52},
\begin{align*}
\wh X(t)&=\wh X_0+\wh A(t) -\sum_{i=0}^2W^{(i)}(t)
+n^{-1/2}\Big[\la^nt-\sum_{i=0}^2T^{(i)}(t)\Big]
\\
&=\wh
X_0+W(t)+n^{-\Half}\Big[\la^nt-\sum_{k=1}^n\mu_k\int_0^tB_k(s)ds\Big],
\end{align*}
where \eqref{31}, \eqref{03} and \eqref{75} are used in the second
equality. Since $B_k=1-I_k$,
$$
 \wh X(t)=\wh X_0+W(t)+n^{-\Half}\Big[\la^n-\sum_{k=1}^n\mu_k\Big]t
 +n^{-\Half}\sum_k\mu_k\int_0^tI_k(s)ds.
$$
By \eqref{21}, $\mu=\la$, hence by \eqref{32} the penultimate term
above is equal to $bt$. This shows \eqref{30}.
 \qed

\skp

\noindent {\bf Acknowledgements:} Adam Shwartz holds The Julius M.\
and Bernice Naiman Chair in Engineering. Work of both authors is
supported in part by the fund for promotion of research and the fund
for promotion of sposored research at the Technion.

\bibliographystyle{plain}

\end{document}